\documentclass[11pt]{amsart}
\usepackage{amsbsy,amsfonts,amsmath,amssymb,eucal,mathrsfs}

\def\op{{\oplus}}
\def\ot{{\otimes}}
\def\ra{{\rightarrow}}
\def\fa{{\mathfrak a}}
\def\fg{{\mathfrak g}}
\def\fh{{\mathfrak h}}
\def\fn{{\mathfrak n}}
\def\fs{{\mathfrak s}}
\def\fsl{{\mathfrak sl}}
\def\PP{{\mathbb P}}
\def\ZZ{{\mathbb Z}}
\def\BC{{\mathbb C}}

\def\cO{{\mathcal O}}
\def\to{\rightarrow}

\newtheorem{prop}{Proposition}

\begin{document}

\title[Homogeneous contact manifold]{Quasi-complete homogeneous contact \nolinebreak manifold \linebreak
associated to a cubic form  }

\author[Jun-Muk Hwang and Laurent Manivel]{Jun-Muk Hwang${}^1$, Laurent Manivel}
\address{Jun-Muk Hwang:  Korea Institute for Advanced Study, Hoegiro 87, Seoul, 130-722, Korea} \email{jmhwang@kias.re.kr}
\address{Laurent Manivel: Institut Fourier, UMR 5582 CNRS/UJF, Universit\'e Joseph Fourier, 38402 Saint Martin d'H\`eres,
France} \email{Laurent.Manivel@ujf-grenoble.fr}
\thanks{${}^1$ supported by the Korea Research Foundation Grant
funded by the Korean Government  (MOEHRD)(KRF-2006-341-C00004)}

 \maketitle

\hspace{3in} {\it Dedicated to S. Ramanan }

\section{Introduction}

This note is at the crossroad of two different lines of study.

On the one hand, we propose a general construction of a
homogeneous quasi-projective manifold $X_c$ associated to a cubic
form with a mild genericity property. These manifolds are
rationally chain connected (Proposition 2), a property which
relates our study to that of certain types of homogeneous spaces
considered in \cite{bb1,bb2,bbk}.

On the other hand, we show that our manifolds $X_c$ are endowed with
natural contact structures (Proposition 3). Our construction thus
appears as part of the general study of contact projective and
quasi-projective manifolds. Of course the projective case is the
most interesting one,  the main open problem in this area being the
Lebrun-Salamon conjecture: the only  Fano contact manifolds should
be the projectivizations of the minimal nilpotent orbits in the
simple Lie algebras. As explained in section 4, our construction is
in fact modeled on these homogeneous contact manifolds, which are
known to be associated to very special cubic forms: the determinants
of the simple cubic Jordan algebras.

Under both points of view, one of the most interesting questions
one may ask about the quasi-projective contact manifolds $X_c$ is
about their compactifications. Even the existence of a small
compactification (that is, with a boundary of codimension at least
two) is not clear. Also, it is extremely tempting to try to
construct new projective contact manifolds by compactifying some
$X_c$ in such a way that the contact structure extends. We show
that this is possible if and only if the cubic $c$ is the
determinant of a simple cubic Jordan algebra (Proposition 6). This
can be interpreted as an evidence for the Lebrun-Salamon
conjecture.

\section{Homogeneous spaces defined from cubics}

Let $V$ be a complex vector space of dimension $p$. Let $c \in S^3
V^*$ be a cubic form on $V$. Let $B : S^2 V \to V^*$ be the system
of quadrics defined by $$B(v_1, v_2) = c( v_1, v_2, \cdot).$$
In all the sequel we make the following

\medskip
\noindent {\sc Assumption on $c$}. {\it The homomorphism $B$ is surjective.}

\medskip
Let $W$ be a complex vector space of dimension 2. Fix a choice of
a non-zero 2-form $\omega \in \wedge^2 W^*$.
%$SL(W)$ denotes the group of
%automorphisms of $W$ preserving $\omega$.

%\medskip
 Let $\fn := \fn_1 \oplus  \fn_2 \oplus  \fn_3$ where
$$ \fn_1 := V \otimes W, \qquad \fn_2 := V^*,
\qquad \fn_3 := W .$$
Define a graded Lie algebra structure on
$\fn$, by
\begin{eqnarray*}
[v_1 \otimes w_1, v_2 \otimes w_2] &= &\omega(w_1, w_2) B(v_1, v_2), \\
 \, [v_1^*, v_2 \otimes w_2 ] &= &v^*_1(v_2) \; w_2.
\end{eqnarray*}
The Jacobi identity holds because $\dim W =2$.

 Let $N$ be the nilpotent Lie group with Lie
algebra $\fn$. For a point $\ell \in \PP W$, denote by $\hat{\ell}
\subset W$ the corresponding 1-dimensional subspace. Let
$$\fa_{\ell} := V \otimes \hat{\ell} \; \subset \; V \otimes W = \fn_1$$
be the abelian subalgebra of $\fn$ and $A_{\ell} \subset N$ be the
corresponding additive  abelian subgroup. We have the smooth
subvariety $A \subset N \times \PP W$ defined by $$A : = \{ (g,
\ell), g \in A_{\ell}\}.$$ This variety $A$ can be viewed as a
family of abelian subgroups parametrized by $\PP W$. Let $\psi:
X_c \to \PP W$ be the family of relative quotients $$X_c := \{
N/A_{\ell}, \ell \in \PP W \}$$ with the quotient map $\xi: N
\times \PP W \to X_c$. Then
$$ \dim X_c = \dim N + 1 - p = 2p +3.$$
Observe that $\psi$ is a locally trivial fibration whose fibers are
isomorphic to affine spaces. But it is not a vector nor an affine bundle.
In fact the transition functions are quadratic, because the nilpotence
index of $N$ is three.

The variety $X_c$ is homogeneous under the action of the group $$G
:=N \triangleleft SL(W) \;\; \mbox{(semi-direct product).}$$  Let
$o \in N$ be the identity and $\ell \in \PP W$ be a fixed base
point. Then $x_\ell := \xi( o \times \ell) $ will be our base
point for $X_c$. Its stabilizer
%of the point of $X_c$ defined by $A_\ell$, for $\ell\in \PP W$,
is $H=A_\ell\triangleleft B_\ell$, if $B_\ell$ denotes the
stabilizer of $\ell$ in $SL(W)$. Moreover a Borel subgroup of $G$
is $B=N\triangleleft B_\ell$, and we have a sequence of quotients
$$G\ra G/B_\ell\stackrel{\xi}{\ra} G/H=X_c\stackrel{\psi}{\ra} G/B=\PP W.$$

\begin{prop}
\begin{enumerate}
\item $X_c$ is simply connected. \item Let $L:= \psi^* {\mathcal
O}_{\PP W}(1)$. Then ${\rm Pic}(X_c) = {\ZZ}L$.
\end{enumerate}
\end{prop}

\proof As a variety, each $A_\ell$ is nothing but an affine space.
So the variety  $X_c$ being
fibered in simply connected manifolds over the projective line, is simply
connected. This proves (1).

The character group $X(G)$ of $G$ being trivial, the forgetful map
$$\alpha : Pic^G(X_c)\ra Pic(X_c)$$ is injective (\cite{mumford}
Proposition 1.4). Moreover the Picard group of $G$ is trivial, so
$\alpha$ is in fact an isomorphism (see the proof of Proposition
1.5 in \cite{mumford}). But $Pic^G(X_c)\simeq X(H)$ and an easy
computation shows that $X(H)=X(B_\ell)$. This implies (2). \qed

\medskip
Note that $G$ is generated by $H$ and $SL(W)$ such that $H \cap
SL(W)$ is a Borel subgroup of $SL(W)$. Thus we can apply
 Proposition 4.1 in \cite{bbk} to deduce:

 \begin{prop}
The variety $X_c$ is rationally chain connected. In particular,
$X_c$ is quasi-complete, i.e.,  there is no non-constant regular
function on $X_c$.
\end{prop}

In fact, it is easy to show that for any $n\in N$, the image of
$\{n\}\times \PP W$ under $\xi$ is a smooth rational curve on
$X_c$ with normal bundle of the form ${\mathcal O}(1)^p \oplus
{\mathcal O}^{p +2}.$

\section{Contact structures}

Consider the tangent spaces
$$T_{o \times \ell}(N \times \PP W) = \fn_1 \oplus  \fn_2
\oplus  \fn_3 \oplus  T_{\ell}(\PP W)$$
$$ T_{x_\ell}(X_c) = \fn_1/\fa_{\ell}
\oplus  \fn_2 \oplus  \fn_3 \oplus  T_{\ell}(\PP W).$$ Using the subspace
$\hat{\ell} \subset W = \fn_3$, we define the hyperplane
$$D_{x_\ell} : = \fn_1/ \fa_{\ell} \oplus  \fn_2 \oplus  \hat{\ell} \oplus  T_{\ell}
(\PP W)$$ inside $T_{x_\ell}(X_c)$. This hyperplane is invariant
under the action of the stabilizer $H$ of $x_\ell$ in $G = N
\triangleleft SL(W)$, so we get a well-defined hyperplane
distribution $D \subset T(X_c)$ with $T(X_c)/D \cong L$.

\begin{prop}
The distribution $D \subset T(X_c)$
defines a contact structure.
\end{prop}

\noindent {\it Remark}. Observe that since
$$ D_{x_\ell} = (V \otimes W/\hat{\ell}) \; \oplus V^* \oplus \hat{\ell}
\oplus {\rm Hom}(\hat{\ell}, W / \hat{\ell}),$$ there is a natural
$W/\hat{\ell}$-valued symplectic pairing on $D_{x_\ell}$. Note
that $W/ \hat{\ell} $ is the fiber of $L$ at $x_\ell$. This shows
that the bundle $D$ has an $L$-valued symplectic form, and indeed
this symplectic form comes from the contact structure
$$0 \to D \to T(X_c) \to L \to 0.$$

\proof Let $Y_c$ be the variety defined as the complement
$L^\times$ of the zero section in the total space of the line
bundle dual to $L$. Let $\theta$ be the $L$-valued 1-form on $X_c$
defining $D$. To check that $\theta$ is a contact form, it
suffices to show that the 2-form $d \tilde{\theta}$ where
$\tilde{\theta}$ is the pull-back of $\theta$ to $Y_c$, is
symplectic  (see \cite{beauville}, Lemma 1.4).

To check this we make a local computation. Let $m\in\PP W$ be some
point distinct from $\ell$ and $\hat{m} \subset W$ be the
corresponding line. Then
$$\fn_m:= (V\ot \hat{m}) \; \op V^*\ot W \; \subset \fn$$ defines a complement
to $\fa_{\ell} \subset \fn$. We can define a local analytic chart
on $X_c$ around $x_\ell$ by
$$x(X,p)=\xi (exp(X)\times p),$$ where $X \in \fn_m$ and $p \in \PP W-m$. Let us write down
 $\theta$ in that local chart. Since the chart
preserves the fibration over $\PP W$ we just need to compute over
$\ell$. The differential $e_X$ of the exponential map at $X$, seen
as an endomorphism of $\fn_m$, is defined by the relation
$$exp(X+te_X(Y)+O(t^2))x_\ell=exp(X)exp(tY)x_\ell.$$ Now we can use
the fact that $N$ being 3-nilpotent, the Campbell-Hausdorff
formula in $N$ is quite simple: we have $exp(X)exp(Y)=exp(H(X,Y))$
for $X,Y\in\fn$, with
$$H(X,Y)=X+Y+\frac{1}{2}[X,Y]+\frac{1}{12}[X,[X,Y]]+\frac{1}{12}[Y,[Y,X]].$$
We easily deduce that $e_X(Y)=Y+\frac{1}{2}[X,Y]+\frac{1}{12}[X,[X,Y]]$.
Now we can decompose this formula with respect to the three-step grading
of $\fn$. If $Z=e_X(Y)=Z_1+Z_2+Z_3$, we find that
\begin{eqnarray*}
Z_1 &= &Y_1, \\
Z_2 &= &Y_2+\frac{1}{2}[X_1,Y_1], \\
Z_3 &= &Y_3+\frac{1}{2}[X_1,Y_2]+\frac{1}{2}[X_2,Y_1]
+\frac{1}{12}[X_1,[X_1,Y_1]],
\end{eqnarray*}
which can be inverted as
\begin{eqnarray*}
Y_1 &= &Z_1, \\
Y_2 &= &Z_2-\frac{1}{2}[X_1,Z_1], \\
Y_3 &= &Z_3-\frac{1}{2}[X_1,Z_2]-\frac{1}{2}[X_2,Z_1]
+\frac{5}{12}[X_1,[X_1,Z_1]].
\end{eqnarray*}
Since the hyperplane $D_{x_\ell}$ is defined by the condition that
$Y_3$ belongs to $\hat{\ell}$, we deduce that the contact form is
given at $x(X,\ell)$, in our specific chart, by the formula
$$\theta_{x(X,\ell)}(Z)=Z_3-\frac{1}{2}[X_1,Z_2]-\frac{1}{2}[X_2,Z_1]
+\frac{5}{12}[X_1,[X_1,Z_1]] \quad \mod\; \hat{\ell}.$$ Even more
explicitly, if we write $Z_1=z_1\ot m$ and $X_1=x_1\ot m$, we have
$[X_1,Z_1]=0$, $[X_1,Z_2]=Z_2(x_1)m$ and $[X_2,Z_1]=-X_2(z_1)m$, so
$$\theta_{x(X,\ell)}(Z)=Z_3+\frac{1}{2}(X_2(z_1)-Z_2(x_1))m.$$
Now we pull-back $\theta$ to $Y_c=L^\times$. A local section of
$L^\times$ around $\ell$
is given by $m^*-z\ell^*$  over the point $p=\ell+zm$ of $\PP W$.
Over $\phi=y(m^*-z\ell^*)$, we get the 1-form on $L^\times$ given in our
local chart by
$$\tilde{\theta}_{x(X,p),\phi}(Z,Y)=y(m^*-z\ell^*)(Z_3)
+\frac{y}{2}(X_2(z_1)-Z_2(x_1)).$$
If $Z_3=Z_3^1m+Z_3^2\ell$, this can also be written as:
$$\tilde{\theta}_{x(X,p),\phi}=y(dX_3^1-zdX_3^2)
+\frac{y}{2}(X_2dX_1-X_1dX_2).$$
We can easily differentiate this expression and evaluate it at $x_\ell$.
We obtain
$$d\tilde{\theta}_{x_\ell,ym^*}=dy\wedge dX_3^1-ydz\wedge dX_3^2
+ydX_2\wedge dX_1.$$
Since $y$ is a non zero scalar this 2-form is everywhere non-degenerate.
By homogeneity this remains true over the whole of $L^\times$,
and the proof is complete. \qed

\section{Projective homogeneous contact varieties}

Consider a  complex simple Lie algebra $\fg$ and the adjoint variety
$$Y_{\fg}=\PP \cO_{min}\subset\PP\fg,$$ the projectivization of the
minimal nilpotent orbit $\cO_{min}$. Then $Y_{\fg}$ is homogeneous
under the action of the adjoint group ${\bf G}=Aut(\fg)$. Suppose
that $Pic(Y_{\fg})\simeq\ZZ$ (this is the case if and only if $\fg$
is not if type A). Then the variety $F$ of lines on $Y_{\fg}$ is
${\bf G}$-homogeneous and we can describe a line as follows. Choose
$T\subset B\subset {\bf G}$ a maximal torus and a Borel sugbroup.
Let $\fg_{\psi}$ denote the root space in $\fg$ associated to the
highest root $\psi$. Then $Y_{\fg}={\bf G} \fg_{\psi}$ and the
stabilizer of $\fg_{\psi}$ is the maximal parabolic subgroup
$P_{\alpha}$ of ${\bf G}$ defined by the unique simple root $\alpha$
such that $\psi-\alpha$ is a root. Moreover the line $\ell =\langle
\fg_{\psi},\fg_{\psi-\alpha}\rangle$ is contained in the adjoint
variety $Y_{\fg}$, and $F={\bf G}.\ell$.

There is a five-step grading on $\fg$ defined by the highest root $\psi$,
as follows. Define $H_\psi\in [\fg_{\psi},\fg_{-\psi}]$ by the condition
that $\psi(H_\psi)=2$. Then the eigenvalues of $ad(H_\psi)$ are
$0,\pm 1,\pm 2$ and the eigenspace decomposition yields the five-step grading
$$\fg=\fg_{-2}\oplus \fg_{-1}\oplus \fg_{0}\oplus \fg_{1}\oplus \fg_{2}.$$
We have $\fg_{2}=\fg_{\psi}$, while $\fg_{\alpha}$ and $\fg_{\psi-\alpha}$
are respectively lines of lowest and highest weights in $\fg_1$.

Since $\fg_{\psi-\alpha}$ defines, exactly as $\fg_{\psi}$, a point of
the adjoint variety, we can use the root $\psi-\alpha$ to define another
five-step grading. Since $ad(H_\psi)$ and $ad(H_{\psi-\alpha})$ commute,
we get a double grading on $\fg$. Moreover,
the stabilizer $\fs\subset\fg$ of the line $\ell \subset\fg$
decomposes as follows (where the grading defined by $ad(H_\psi)$
can be read horizontally):
$$\begin{array}{cccc}
\fg_{-\alpha} & & & \\
\fg_{00} & \fg_{10} & & \\
\fg_{\alpha} & \fg_{11} & \fg_{21} & \fg_{\psi-\alpha} \\
 & & & \fg_{\psi}
\end{array}$$
Let $W=\fg_{\psi-\alpha}\oplus\fg_{\psi}\simeq \BC^2$ and
$V=\fg_{21}^*$. The map $\fg_{10}\otimes \fg_{21}\rightarrow
\fg_{31}=\fg_{\psi-\alpha}$ defined by the Lie bracket is a
perfect pairing, as well as $\fg_{11}\otimes \fg_{21}\rightarrow
\fg_{32}=\fg_{\psi}$, giving  a natural identification
$$\fg_{10}\oplus\fg_{11}\simeq V\ot W $$ and isomorphisms $$ \fg_{10}
\stackrel{\phi}{\longrightarrow} \fg_{11} \cong V.$$ The positive
part of the vertical grading of $\fs$ thus reads $$(V\ot W) \; \op
V^*\op W=\fn.$$ Note that the degree zero part of this grading
reads $\fsl(W)\times \fh_{00}$, where $\fg_{00}=[\fg_{-\alpha},
\fg_{\alpha}]\op \fh_{00}$ is an orthogonal decomposition with
respect to the Killing form.

The cubic form $c$ on $V$ is defined (up to scalar) once we
identify $\fg_{10}$ with $\fg_{11}$, through the map $\phi$. We
also need to choose a generator $X_{\psi}$ of $\fg_{\psi}$.   Then
we can define $c$ by the formula
$$[\phi(X), [\phi(X), X]]=c(X)X_{\psi} \qquad\forall X\in \fg_{10}\simeq V.$$

\noindent {\it Remark}. This construction is closely  related to the
{\it ternary models} for simple Lie algebras considered in \cite{man},
section 2. These models are of the form
$$\fg = \fh\times\fsl(U)\; \op\;  (U\ot V)\;  \op\;  (U^*\ot V^*),$$
where $U$ is three dimensional, and $V$ is an $\fh$-module.
To define a Lie bracket on $\fg$, one needs a cubic form $c$ on $V$,
a cubic form $c^*$ on $V^*$, and a map $\theta : V\ot V^*\ra\fh$.
Then the Jacobi identity implies a series of conditions on these
data, including that
$$\fh\subset Aut(c)\cap Aut(c^*).$$
These conditions should ultimately lead to a cubic
Jordan algebra structure on $V$. If we choose a maximal torus in
$\fsl(U)$ and use the associated grading on $U,U^*$, we get an
hexagonal model as in Figure 2 of \cite{mukai} :

$$\begin{array}{ccccccc}
 & & & \fg_{-\alpha} & & & \\ &&&&&&\\
\fg_{-\psi} & & \fg_{-1-1} & & \fg_{10} & & \fg_{\psi-\alpha}\\&&&&&&\\
 & \fg_{-2-1} & & \fg_{00} & & \fg_{21} & \\&&&&&&\\
\fg_{\alpha-\psi} & & \fg_{-10} & & \fg_{11} & & \fg_{\psi}\\&&&&&&\\
 & & & \fg_{\alpha} & & &
\end{array}$$

The subalgebra we denoted $\fn$ is the sum of the factors in the last
three columns. and we can add the factor $\fsl(W)$
from the middle column in order to get $\fs$.

\medskip Once we have defined the cubic $c$ associated to the
simple Lie algebra $\fg$, we have the associated homogeneous space
$X_c$ with its natural contact structure. A direct verification
gives:

\begin{prop}
The homogeneous space $X_c$ is  an open subset of $Y_\fg$, with a
codimension two boundary.
%swept out
%by the translates of the line $\ell$ under the parabolic group $P$.
Its contact structure is the restriction of the natural contact
structure on $Y_\fg$.
\end{prop}

It is tempting, but illusory, as we shall see, to try to construct
new projective contact manifolds as suitable compactifications of our
homogeneous spaces $X_c$ for other types of  cubics.

\section{Compactifications}

Since all regular functions on $X_c$ are constant, we can expect
that $X_c$ admits a {\it small compactification}, that is, a
projective variety $\bar X_c$ containing $X_c$ as an open subset
in which the boundary of $X_c$ has codimension two or more. By
Theorem 1 in \cite{bb2} and Lemma 1, it is enough to check that
the algebra of sections
$$R(X_c,L)=\bigoplus_{k=0}^\infty \Gamma(X_c,L^k)$$
is of finite type, as well as all the $R(X_c,L^m)$, for $m\ge 1$.
We have not been able to prove this but we can make the
following observations.

\medskip
Since the Lie algebra $\fg$ of $G = N \triangleleft SL(W)$
preserves the contact structure we have defined on $X_c$,
there must be a morphism
$\varphi$ from $X_c$ into $\PP \fg^*$ (see \cite{beauville},
Section 1). In fact, any contact vector
field on a contact manifold defines a holomorphic section of the
contact line bundle $L$. Thus $\fg$ defines a linear subsystem in
$|L|$. The morphism is always etale over its image, and since our
$X_c$ is simply connected we conclude that $\varphi$ embeds $X_c$
as a coadjoint orbit in $\fg^*$. We thus have a natural projective
compactification of $X_c$ in $\PP \fg^*.$

Note that the inclusion of  $X_c$ in $\PP \fg^*$ is just the
projectivization of the moment map of the symplectic variety
$Y_c=L^\times$. We have a commutative diagram
$$\begin{array}{ccc}
 Y_c & \stackrel{\mu}{\longrightarrow} & \fg^* \\
\downarrow  & & \downarrow \\
 X_c &  \stackrel{\nu}{\longrightarrow} & \PP\fg^*
\end{array}$$
 \medskip
Here $\mu$ denotes  the ($G$-equivariant) moment map and $\nu$ is
its quotient by the ${\BC}^*$-action. We have $\fg^*=
\fn^*\;\op\;\fsl(W)^*$ and the component $\mu'$ of $\mu$ on
$\fn^*$ is not injective, since the $N$ action on $X_c$ preserve
the $\PP^1$-fibration. Consider $\mu'(Y_c)\subset\fn^*$.

\begin{prop}
Suppose that the cubic hypersurface $Z_c\subset\PP V$ be smooth.
Then the boundary of $\mu'(Y_c)$ has codimension at least two.
\end{prop}

\proof We can describe explicitly the closure of $\mu'(Y_c)$  as
the set of triples $(\phi_1,\phi_2,\phi_3) \in\fn^*$ such that
$$\omega (\phi_1,\phi_3)=c(\phi_2,\phi_2,.),$$
where $\omega :(V^*\ot W)\times W\rightarrow V^*$ is the natural
bilinear map.

If $\phi_3\ne 0$, we are in $\mu'(Y_c)$. Thus on the boundary, we must
have $\phi_3=0$, and then $c(\phi_2,\phi_2,.)=0$. But under our
smoothness assumption on $Z_c$, this implies that $\phi_2=0$.
So the boundary of $\mu'(Y_c)$ has dimension at most $2p=\dim\mu'(Y_c)-2$,
the number of parameters for $\phi_1$.
\qed

\medskip
This seems to be a first step towards proving that $X_c$ has a small
compactification. But we have not been able even to find conditions on $c$
that would ensure that the compactification $\bar X_c\subset\PP\fg^*$
is small.

\smallskip
What is rather surprising is that the cubics whose associated
variety $X_c$ has a smooth contact compactification can be
completely classified. By this, we mean a smooth projective
variety $\bar X_c$ compactifying $X_c$, with a contact structure
extending that of $X_c$.

\begin{prop}
There exists a smooth contact compactification $\bar{X}_c$ of
$X_c$  if and only if $c$ is the cubic norm of a semi-simple
Jordan algebra.
\end{prop}

\proof
We will deduce this statement from a study the variety of minimal rational
tangents ${\mathcal C}_{x_\ell} \subset \PP D_{x_\ell}$. Note that the
space of lines on $X_c$ through $x_\ell$ is just $\xi^{-1}(x_\ell) \cong
A_{\ell} \cong V$.

\medskip
We claim that the tangent map sending a line through $x_\ell$ to its
tangent direction in $\PP D_{x_\ell}$ is equal, up to scalars,
 to the rational
map $$ \tau: V \to \PP D_{x_\ell} = \PP (  V \oplus V^* \oplus
{\BC} \oplus {\BC})$$ $$ \tau(v) := \; [ v: B(v,v): c(v,v,v):1].$$
Indeed, a line through $x_\ell$ in $X_c$ is of the form
$\ell_g=\xi (g\times\PP W)$ for $g\in A_{\ell}$. To write this
line in the local chart we used in the proof of Proposition 3 (we
use the same notations), we must write
$$\xi (g\times p)=exp(Z)\xi (o\times p).$$
If $g=exp(X)$ with $X\in V\ot \ell$, this amounts to solving the equation
$exp(X)=exp(Z)exp(W)$, with $Z\in \fn_m$ and $W\in V\ot p$.
So $X=H(Z,W)$, and if we write $X=v\ot \ell$ for some $v\in V$,
we must have $W=v\ot p$ and then we get, up to term of order at least two
in $z$,
\begin{eqnarray*}
Z_1 &= &-zv\ot m, \\
Z_2 &= &-\frac{1}{2}[Z_1,W]=\frac{z}{2}B(v,v), \\
Z_3 &= &-\frac{1}{6}[Z_1,[Z_1,W]]=-\frac{z}{6}c(v)\ell.
\end{eqnarray*}
This proves the claim.

We can now conclude the proof as follows. By the results of
\cite{keb}, the closure of the image of this map must be smooth if
there exists a smooth contact compactification of $X_c$. So the
closure of the image of the map $\tau$ must be smooth. But then we
can apply Corollary 26 in  \cite{lm}. \qed


\begin{thebibliography}{aa}

\bibitem{beauville}
Beauville A.:
{\it Fano contact manifolds and nilpotent orbits},
Comment. Math. Helv. {\bf 73} (1998), no. 4, 566--583.

\bibitem{bb1}
Bien F., Borel A.: {\it Sous-groupes \'epimorphiques de groupes lin\'eaires
alg\'ebriques} I,  C. R. Acad. Sci. Paris S\'er. I Math. {\bf 315}  (1992),
no. 6, 649--653.

\bibitem{bb2}
Bien F., Borel A.: {\it Sous-groupes \'epimorphiques de groupes lin\'eaires
alg\'ebriques} II,
C. R. Acad. Sci. Paris S\'er. I Math. {\bf  315}  (1992),  no. 13, 1341--1346.

\bibitem{bbk}
Bien F., Borel A., Koll\'ar J.: {\it
Rationally connected homogeneous spaces},
Invent. Math. {\bf 124} (1996), no. 1-3, 103--127.

\bibitem{keb}
Kebekus S.: {\it
Lines on complex contact manifolds} II,
Compos. Math. {\bf 141} (2005), no. 1, 227--252.

%\bibitem{ }[ ]

\bibitem{lm}
Landsberg J.M., Manivel L.: {\it Legendrian varieties},  Asian J. Math.
{\bf 11}  (2007),  no. 3, 341--359.

\bibitem{man}
Manivel L.: {\it Configurations of lines and models of Lie algebras},
J. Algebra {\bf 304}  (2006),  no. 1, 457--486.

\bibitem{mukai}
Mukai S.: {\it Simple Lie algebra and Legendre variety}, preprint 1998.

\bibitem{mumford} Mumford D., Fogarty J., Kirwan F.:
Geometric invariant theory, Third edition, Springer 1994.


%\bibitem{ }[ ]


\end{thebibliography}
\end{document}